\documentclass[11pt]{article}
\usepackage[dvips]{graphicx}
\usepackage{amssymb}
\usepackage{amsmath}
\usepackage{amsthm}
\usepackage{latexsym}
\usepackage[english]{babel}

\newtheorem{defi}{Definition}[section]
\newtheorem{teo}[defi]{Theorem}

\newtheorem{lem}[defi]{Lemma}

\newtheorem{es}[defi]{Example}

\begin{document}


\title{A symmetry result on Reinhardt domains}
\author{Vittorio Martino$^{(1)}$ } \addtocounter{footnote}{1}
\footnotetext{Dipartimento di Matematica, Universit\`a di Bologna,
piazza di Porta S.Donato 5, 40127 Bologna, Italy. E-mail address:
{\tt{martino@dm.unibo.it}}}
\date{}
\maketitle
$$$$
{\noindent\bf Abstract} We show the
following symmetry property of a bounded Reinhardt domain $\Omega$ in $\mathbb{C}^{n+1}$:
let $M=\partial\Omega$ be the smooth boundary of $\Omega$ and let $h$ be the Second Fundamental Form of $M$;
if the coefficient $h(T,T)$ related to the characteristic direction $T$ is constant then $M$ is a sphere.
In Appendix we state the result from an hamiltonian point of view.
$$$$
$$$$
\section{Introduction}
A Reinhardt domain $\Omega$ (with center at the origin) is by definition an open subset of $\mathbb{C}^{n+1}$ such that
\begin{equation}\label{reinhardtdefinition}
\mbox{if} \quad (z_1,\ldots,z_{n+1})\in \Omega \quad \mbox{then} \quad (e^{i\theta_1} z_1,\ldots, e^{i \theta_{n+1}} z_{n+1} )\in \Omega
\end{equation}
for all the real numbers $\theta_1,\ldots,\theta_{n+1}$.
These domains naturally arise in the theory of several complex variables as the logarithmically convex Reinhardt domains are the domains of convergence of power series (see for instance \cite{hormandercomplex}, \cite{jarpfl}). We will suppose from now on that the Reinhardt domain $\Omega$ has a smooth boundary (it would be enough $C^2$).  The boundary $M:=\partial\Omega$ is then a smooth real hypersurface in $\mathbb{C}^{n+1}$ and thus a CR-manifold of CR-codimension equal to one, with the standard CR structure induced by the holomorphic structure of $\mathbb{C}^{n+1}$. Thus for every $p\in M$ the tangent space $T_p M$ splits in two subspaces: the $2n-$dimensional horizontal subspace $H_p M$, the largest subspace in $T_p M$ invariant under the action of the standard complex structure $J$ of $\mathbb{C}^{n+1}$ and the vertical one-dimensional subspace generated by the characteristic direction $T_p :=J \cdot N_p$, where $N_p$ is the unit normal at $p$. Moreover, if $\widetilde{g}$ is the standard metric on
$\mathbb{C}^{n+1}$, then it holds
$$T_p M=H_p M\oplus\mathbb{R}T_p$$
and the sum is $\widetilde{g}$-orthogonal.\\
Let us consider the complexified horizontal space
$$H^{\mathbb{C}} M :=\{Z= X-iJ \cdot X: X\in HM \}$$
The Levi Form $l$ is then the sesquilinear and hermitian operator on $H^{\mathbb{C}}M$ defined in the following way:
$\forall Z_1,Z_2\in H^{\mathbb{C}} M$
\begin{equation}\label{leviform}
l(Z_1,Z_2)= \widetilde{g}({\widetilde{\nabla}}_{Z_1} \bar Z_2, N)
\end{equation}
where $\widetilde{\nabla}$ is the Levi-Civita connection for $\widetilde{g}$.
Moreover by a direct computation it holds
\begin{equation}\label{leviformquadratic}
l(Z,Z)= \widetilde{g}({\widetilde{\nabla}}_{Z} \bar Z, N)=\widetilde{g}([X,Y],T)
\end{equation}
where $Y=J \cdot X$. We will say $M$ be (strictly) pseudoconvex if $l$ is (strictly) positive definite as quadratic form.\\
In analogy with classical curvatures defined in terms of elementary symmetric functions of the eigenvalues of the Second Fundamental Form, one defines the $j$-th Levi curvatures  $L^{j}$ in terms of elementary symmetric functions of the eigenvalues of the Levi Form
$$L^{j}= \frac{1}{\binom {n}{j}} \sum_{1\leq i_1<\dots<i_j\leq n}\lambda_{i_1}\cdots\lambda_{i_j},$$
where  $\lambda_1,\dots, \lambda_{n}$ are the eigenvalues of $l$. In particular when $j=n$ we have the Total-Levi Curvature and when $j=1$ we have the Levi-Mean Curvature $L$.\\
Being hypersurfaces in $\mathbb{C}^{n+1}$ real hypersurfaces in $\mathbb{R}^{2n+2}$, one can also compare the Levi Form with the Second Fundamental Form $h$ of $M$ by using the identity \cite{bog}
$$
l(Z,Z)=h(X,X)+h(J(X),J(X)),\quad \forall X\in HM
$$
Thus, a direct calculation leads to the relation between the classical Mean Curvature $H$ and the Levi-Mean Curvature $L$ \cite{io}:
\begin{equation}\label{meanlevimean}
H=\frac{1}{2n+1}(2nL+h(T,T))
\end{equation}
where $h(T,T)=\widetilde{g}({\widetilde{\nabla}}_{T} T, N)$ is the coefficient of the Second Fundamental Form related to the characteristic direction $T$.
\begin{defi}\label{linea}
We will call $h(T,T)$ the characteristic curvature of $M$.
\end{defi}

\noindent
By (\ref{meanlevimean}) the characteristic curvature is a sort of complementary of the Levi-Mean Curvature in computing the Mean Curvature. Moreover, for every hypersurface in $\mathbb{C}^{n+1}$, $h(T,T)$ is invariant under a biholomorphic (rigid) transformation, as the Levi curvatures are.\\
Following the pioneering result due to Alexandrov \cite{alex} on the classical Mean Curvature of Euclidean surface, the problem of characterizing compact hypersurfaces with positive constant Levi-Mean Curvature has recently received a great amount
of attention. Klingenberg in \cite{kling} gave a first positive answer to this problem by showing that if the
characteristic direction is a geodesic and the Levi Form is diagonal, then $M$ is a sphere. Monti and Morbidelli in \cite{momo} proved a Darboux-type theorem for $n\geq 2$: the unique Levi umbilical hypersurfaces in $\mathbb{C}^{n+1}$ with all constant Levi curvatures are spheres or cylinders. Later on Montanari and the author proved two results of this type: in \cite{MaM} they relaxed Klingerberg conditions and they proved that if  the characteristic direction is a geodesic, then Alexandrov Theorem holds for hypersurfaces with positive constant Levi-Mean Curvature; in \cite{marmon} they proved some integral formulas for compact hypersurfaces, of independent interest, and then they follow the Reilly approach \cite{RE2}, \cite{R}, \cite{R2} to prove Isoperimetric estimates and a Alexandrov type theorem, namely: let $M$ be a closed smooth real hypersurface bounding a star-shaped domain in $\mathbb{C}^{n+1}$, if the $j$-Levi curvature is a positive constant $K$ and the maximum of the Mean Curvature of $M$ is bounded from above by $K$ then $M$ is a sphere.\\
In a couple of recent papers Hounie and Lanconelli proved Alexandrov type theorems for Reinhardt domains in $\mathbb{C}^{2}$ first and for Reinhardt domain in $\mathbb{C}^{n+1}$, $n\geq 1$, with an additional rotational symmetry then. In \cite{HL} they showed the result for bounded Reinhardt domain of $\mathbb{C}^2,$
i.e. for domains $\Omega$ such that if $ (z_1,z_2)\in \Omega$ then $ (e^{i
\theta_1} z_1, e^{i \theta_2} z_2 )\in \Omega$ for all real $\theta_1,
\theta_2 .$ Under this hypothesis, in a neighborhood of a point,
there is a defining function $F$ only depending on the radii
$r_1=|z_1|,$ $r_2=|z_2|,$ $F(r_1,r_2)=f(r_2^2)-r_1^2$ with $f$ the
solution of the ODE
\begin{equation}\label{ODE}
  sff''=sf'^2-k(f+sf'^2)^{3/2}-ff'
\end{equation}
Alexandrov Theorem follows from uniqueness of the solution
of \eqref{ODE}. Their technique has then been used in \cite{HL1}
to prove an Alexandrov Theorem for bounded Reinhardt domains in $\mathbb{C}^{n+1}$
with an additional rotational symmetry in two complementary sets of variables, for every $n$.\\
Here we prove a similar result of symmetry for Reinhardt domains in $\mathbb{C}^{n+1}$ starting from the characteristic curvature rather than the Levi ones.
\begin{teo}\label{teoreinhardt}
Let $M:=\partial\Omega$ be the smooth boundary of a bounded Reinhardt domain $\Omega$ in $\mathbb{C}^{n+1}$.
If the characteristic curvature $h(T,T)$ is constant then $M$ is a sphere of radius equal to $1/h(T,T)$.
\end{teo}

\noindent
Let $\{X_1,\ldots,X_{n},Y_1,\ldots,Y_{n}\}$, with $Y_k=J \cdot X_k$, be an orthonormal basis of the horizontal space $HM$; keeping in mind the structure of the Second Fundamental Form
$$h=\left(
\begin{array}{cccccc}
 h(X_k,X_k) & h(X_k,Y_j)  & h(X_k,T)  \\
 h(Y_j,X_k) & h(Y_j,Y_j)  & h(Y_j,T)  \\
 h(T,X_k)   & h(T,Y_k)    & h(T,T)    \\
\end{array}\right)$$
with $k$ and $j$ running in $1,\ldots,n$, we are making assumption only on the one-dimensional characteristic subspace of the tangent space rather than on the $2n-$dimensional horizontal one $HM$: moreover when in addition one assumes one of the Levi curvatures be non zero (as in the Alexandrov type results) then $HM$ spans the whole tangent space; in fact the vector fields $\{X_1,\ldots,X_{n},Y_1,\ldots,Y_{n}\}$ satisfy the H\"{o}rmander rank condition.\\
When there exists a defining function $f:\mathbb{C}^{n+1}\rightarrow \mathbb{R}$
$$\Omega=\{z\in\mathbb{C}^{n+1}:f(z)<0\},\quad M=\partial\Omega=\{z\in\mathbb{C}^{n+1}:f(z)=0\}$$
such that $f(z)=g(r)$ depends only on the radii $r=(r_1,\ldots,r_{n+1})$, where
$$r_k=z_k \bar z_k, \qquad k=1,\ldots,n+1$$
then we can find an explicit formula to compute the characteristic curvature $h(T,T)$. In fact by using the following identities
$$f_k=\bar z_k g_k, \qquad f_{\bar k}= z_k g_k ,
\qquad f_{\bar j k}=\delta_{j k}g_k+z_j \bar z_k g_{j k}$$
$$|\partial f|^2=\sum_k r_k g_k^2$$
the unit normal $N$ is
$$N=-\frac{1}{|\partial f|} \sum_k (z_k g_k \partial_{z_k}+\bar z_k g_k\partial_{\bar z_k})$$
and the characteristic direction $T$ reads as
$$T=J \cdot N=-\frac{i}{|\partial f|} \sum_k (z_k g_k\partial_{z_k}-\bar z_k g_k\partial_{\bar z_k})$$
Then by a direct computation we have that
\begin{equation}\label{reinhardtcharacteristiccurvature}
h(T,T)=\widetilde{g}({\widetilde{\nabla}}_{T} T, N)=\sum_k^{n+1} \frac{r_k g_k^3}{|\partial f|^{3}}
\end{equation}
\begin{es}[characteristic curvature of the sphere]\label{characteristiccurvaturesphere}
Let
$$g(r_1,\ldots,r_{n+1})=r_1+ \ldots +r_{n+1}-R^2$$
be the defining function of the sphere of radius equal to $R$ in $\mathbb{C}^{n+1}$. By the formula (\ref{reinhardtcharacteristiccurvature}) we have that the characteristic curvature of the sphere is
$$h(T,T)=\frac{1}{R}$$
\end{es}
\begin{es}[characteristic curvature of ellipsoidal type domains]\label{characteristiccurvatureellipsoid}
Let
$$g(r_1,\ldots,r_{n+1})=\frac{r_1}{a_1^2}+ \ldots +\frac{r_{n+1}}{a_{n+1}^2}-1$$
be the defining function of an ellipsoid in $\mathbb{C}^{n+1}$ with $(a_1,\ldots,a_{n+1})$ positive constants. By the formula (\ref{reinhardtcharacteristiccurvature}) we have that at a point $p=(r_1,\ldots,r_{n+1})\in M$ its characteristic curvature  is
$$h_p(T,T)= \frac{\displaystyle \sum_k^{n+1}\frac{r_k}{a_k^6}}{\displaystyle \Big( \sum_k^{n+1}\frac{r_k}{a_k^4} \Big)^{3/2}}$$
\end{es}
\noindent
In the next section we will prove the Theorem \ref{teoreinhardt}, then in the Appendix we will show an Hamiltonian point of view of the result.

\section{Proof of Theorem \ref{teoreinhardt}}
Let us identify $\mathbb{R}^{n+1} \times \mathbb{R}^{n+1} \simeq \mathbb{C}^{n+1}$ so that $z=(x,y)$. First we prove a property of independent interest.
\begin{lem}\label{scalarproductzero}
Let $\Omega$ be a Reinhardt domain in $\mathbb{C}^{n+1}$ and let
$$p=(z_1,\ldots,z_{n+1})=(x_1,\ldots,x_{n+1},y_1,\ldots,y_{n+1})$$
the ``position vector'' of a point on $M:=\partial\Omega$.
If $T_p$ is the characteristic direction at $p \in M$ then it holds identically
\begin{equation}\label{scalprodzero}
\tilde g(p,T_p)\equiv 0
\end{equation}
\end{lem}
\begin{proof}
If $M$ is any smooth hypersurface bounding a domain $\Omega$ in $\mathbb{C}^{n+1}$ with defining function $f:\mathbb{C}^{n+1}\rightarrow \mathbb{R}$ such that
$$\Omega=\{z\in\mathbb{C}^{n+1}:f(z)<0\},\quad M=\partial\Omega=\{z\in\mathbb{C}^{n+1}:f(z)=0\}$$
then the unit normal $N$ is:
$$N=-\frac{1}{|\partial f|}\sum_{k=1}^{n+1}(f_{\bar k}\partial_{z_k}+f_k\partial_{\bar z_k})$$
where $f_k=\displaystyle\frac{\partial f}{\partial z_k}$, with $k=1,\ldots,n+1$.
Thus the characteristic direction $T$ is:
$$T=J \cdot N=-\frac{i}{|\partial f|}\sum_{k=1}^{n+1}(f_{\bar k}\partial_{z_k}-f_k\partial_{\bar z_k})$$
By identifying $f(z)=f(x,y)$, from the real point of view we have:
$$N=-\frac{1}{|\nabla f|}\sum_{k=1}^{n+1}(f_{x_k}\partial_{x_k}+f_{y_k}\partial_{y_k})$$
$$T=\frac{1}{|\nabla f|}\sum_{k=1}^{n+1}(f_{y_k}\partial_{x_k}-f_{x_k}\partial_{y_k})$$
Now, if $\Omega$ is a Reinhardt domain (with center at the origin) in $\mathbb{C}^{n+1}$ then we can find (at least locally) a defining function $f(z)=g(r)$ depending only on the radii $r=(r_1,\ldots,r_{n+1})$ where
$$r_k=z_k \bar z_k = x_k^2 + y_k^2, \qquad k=1,\ldots,n+1$$
So if $g_k=\displaystyle\frac{\partial g}{\partial r_k}$ we obtain
$$f_{x_k}=2 x_k g_k, \qquad f_{y_k}=2 y_k g_k$$
with $k=1,\ldots,n+1$. In vectorial notation then we have
$$T=\displaystyle\frac{1}{|\nabla f|}(f_{y_1},\ldots,f_{y_{n+1}},-f_{x_1},\ldots,-f_{x_{n+1}})=$$
$$=\displaystyle\frac{2}{|\nabla f|}(y_1 g_1,\ldots,y_{n+1} g_{n+1},-x_1 g_1,\ldots,-x_{n+1} g_{n+1})$$
and thus it holds identically
$$\tilde g(p,T_p)=\displaystyle\frac{2}{|\nabla f(p)|} \sum_{k=1}^{n+1}\Big(x_k y_k g_k(p)- y_k x_k g_k(p)\Big)  \equiv 0$$
for every $p\in M$
\end{proof}

\noindent
In other words, the vector position $p$ has generally a normal component and a tangential component; in turn, the tangential component has an horizontal component and a characteristic component: for Reinhardt domains the characteristic component of the vector position $p$ identically vanishes.\\
Now we can prove the main result.
\begin{proof}(of Theorem \ref{teoreinhardt})
Let us consider the function:
$$\varphi:M\rightarrow \mathbb{R}, \qquad \varphi(p)= \frac{|p|}{2}^2=\frac{\tilde g(p,p)}{2}$$
that represents one half the squared distance of the manifold from the origin.
If $V \in TM$ is a tangent vector field to $M$ then the derivative of $\varphi$ along $V$ is
$$V(\varphi(p))=\frac{1}{2}V(\tilde g(p,p))=\tilde g(p,V_p)$$
and by Lemma \ref{scalarproductzero} we have
$$T(\varphi)=\tilde g(p,T)\equiv 0$$
Thus, if $\widehat{p}$ is a critical value of $\varphi$, then
$$X_k(\varphi)_{\mid_{\widehat{p}}}=Y_k(\varphi)_{\mid_{\widehat{p}}}=0$$
Moreover, $\varphi$ evaluated at a critical value is
\begin{equation}\label{varphicriticalvalue}
\varphi(\widehat{p})=\frac{|\widehat{p}|}{2}^2
\end{equation}
and the position vector of any critical value $\widehat{p}$ is parallel to the (inner) unit normal
direction $N$ at $\widehat{p}$
$$\widehat{p}=\tilde g(\widehat{p},N_{\widehat{p}})N_{\widehat{p}}=-|\widehat{p}|N_{\widehat{p}} $$
Differentiating again $\varphi$ along the characteristic direction $T$ we obtain
$$0\equiv T^2(\varphi)=T(\tilde g(p,T))=
\tilde g(T,T)+\tilde g(p,\widetilde{\nabla}_T T)=1+\tilde g(p,\widetilde{\nabla}_T T)$$
and if $\widehat{p}$ is a critical value for $\varphi$ then we get
\begin{equation}\label{criticalvalue}
1-|\widehat{p}| \tilde g(N_{\widehat{p}},\widetilde{\nabla}_T T)=1-|\widehat{p}|h_{\widehat{p}}(T,T)=0
\end{equation}
where $h_{\widehat{p}}(T,T)$ is the characteristic curvature of $M$ at $\widehat{p}$.\\
Since $M$ is a smooth compact hypersurface, then $\varphi$ admits maximum and minimum which are critical values for $\varphi$.
If $h(T,T)$ is constant then by (\ref{criticalvalue}) we have
$$|\widehat{p}|=\frac{1}{h_{\widehat{p}}(T,T)}=\frac{1}{h(T,T)}=const.$$
Then by (\ref{varphicriticalvalue}) $\varphi$ is constant on $M$ and it holds
$$(2\varphi(p))^{1/2}=|p|=\frac{1}{h(T,T)}=const.$$
for every $p \in M$, and it means that $M$ is a sphere of radius $\displaystyle\frac{1}{h(T,T)}$
\end{proof}
\noindent
The boundedness hypothesis is crucial as the next example shows.
\begin{es}[characteristic curvature of a cylinder type domain]\label{characteristiccurvaturecylinder}
Let
$$g(r_1,r_2)=r_1-R^2$$
be the defining function of a cylinder type domain in $\mathbb{C}^{2}$. By the formula (\ref{reinhardtcharacteristiccurvature}) we have that the its characteristic curvature is constant:
$$h(T,T)=\frac{1}{R}$$
\end{es}

\section{Appendix}
Here we want to look at the Reinhardt domains from an hamiltonian point of view. First we recall that for every hypersurface $M$ in $\mathbb{C}^{n+1}$, with $f$ as defining function, the characteristic direction of $M$ is exactly the (normalized) hamiltonian vector field for the hamiltonian function $f$. In fact let us consider a dynamic system with hamiltonian function (smooth enough) depending on position and momentum variables
$$H:\mathbb{R}^{n+1} \times \mathbb{R}^{n+1} \rightarrow \mathbb{R},\qquad z=(q,p)\mapsto H(q,p)$$
and define the Action functional
$$A(z)=\int_{t_0}^{t_1}\Big(\langle p, \dot{q}\rangle -H(q,p)\Big)dt, \qquad z:[t_0,t_1]\rightarrow\mathbb{R}^{2n+2}$$
The first variation of $A$ on a suitable space of curves leads to the following system of differential equations (Hamilton)
\begin{equation}\label{hamilton}
\left\{
\begin{array}{l}
\dot{q}_k=\displaystyle\frac{\partial H}{\partial p_k}(q,p)\\
\\
\dot{p}_k=\displaystyle -\frac{\partial H}{\partial q_k}(q,p)\\
\end{array}\right.
\qquad k=1,\ldots,n+1
\end{equation}
Now, a Least Action Principle states that trajectories of motion (in the generalized phase space $\mathbb{R}^{n+1} \times \mathbb{R}^{n+1}$) are solutions of (\ref{hamilton}). The isoenergetic surface of  $H$ of energy $E$ is the following hypersurface in $\mathbb{R}^{2n+2}$: $M=\{ z \in \mathbb{R}^{2n+2} : H(z)=E \}$. The conservation of energy principle ensures that if $z$ is a critical point for $A$, then $z(t) \in M, \forall t\in [t_0,t_1]$.
The hamiltonian vector field for $H$ is the tangent vector field to $M$
$$X^H_z:=\Big(\frac{\partial H}{\partial p}(q,p),-\frac{\partial H}{\partial q}(q,p)\Big)=J\cdot \nabla H(q,p)$$
where
$$J=\left(
\begin{array}{cc}
0 & I_{n+1}\\
-I_{n+1} & 0\\
\end{array}
\right)
$$
is the canonical symplectic matrix in $\mathbb{R}^{2n+2}$ and in our case it coincides with the standard complex structure in $\mathbb{C}^{n+1}$.\\
The Hamilton system (\ref{hamilton}) rewrites as
$$\dot{z}=X^H_z$$
Now, if one identifies
$$\mathbb{C}^{n+1}\approx\mathbb{R}^{2n+2},\qquad z=(z_1,\ldots,z_{n+1}),\qquad z_k=x+iy\simeq (x_k,y_k)$$
then the hypersurface $M$ defined by
$$M=\{ z \in \mathbb{C}^{n+1} : f(z)=0 \},\qquad f:\mathbb{C}^{n+1} \rightarrow \mathbb{R}$$
is exactly the isoenergetic surfaces of  $H=f+E$. Thus the hamiltonian vector field on $M$ is
$$X^H_z=J\cdot \nabla H(z)= J\cdot \nabla f(z)=J\cdot N=T$$
where $N=\nabla f$ is the normal direction to $M$ and $T$ is the (not normalized) characteristic direction.
Moreover the integral curves of $X^H$ (the orbits in the phase space) coincide with that ones of $T$, eventually reparametrized. In this situation the characteristic curvature $h(T,T)$ is the normal curvature of the hamiltonian trajectories on the isoenergetic surface in the generalized phase space $\mathbb{R}^{n+1} \times \mathbb{R}^{n+1}$.\\
Now, we recall that if $\Omega$ is a Reinhardt domain (with center at the origin) in $\mathbb{C}^{n+1}$ then we can find (at least locally) a defining function $f(z)=g(r)$ depending only on the radii $r=(r_1,\ldots,r_{n+1})$ where
$$r_k=z_k \bar z_k = x_k^2 + y_k^2, \qquad k=1,\ldots,n+1$$
This means that the hamiltonian function depends only on the quantities $r_k=q_k^2+p_k^2$ that represent the actions in the pair of variables action-angle. Thus the angle variables are cyclic and then the actions $r_k$ (and all the functions depending on them) are conserved quantities along the trajectories of motion. In fact we have that the characteristic direction $T$ is:
$$T=-\frac{i}{|\partial f|} \sum_k (z_k g_k\partial_{z_k}-\bar z_k g_k\partial_{\bar z_k})$$
then it holds
$$T(r_k)=0,  \qquad k=1,\ldots,n+1$$
Moreover the system (\ref{hamilton}) reads as
\begin{equation}\label{eq:hamrein}
\dot{z}_k=-if_{\bar k}=-i z_k g_k
\end{equation}
and since $g_k(t)=g_k(0)$, then the curve
$$ z(t)=z_k(0)e^{\displaystyle -ig_k(0)t}$$
is an explicit solution of (\ref{eq:hamrein}) with initial condition $z_k(0)$.\\
In particular, we have that the following curves
$$ z(t)=z_k(0) e^{\displaystyle -i  \frac{ g_k(0)}{|\partial f(0)|}t}$$
are integral curves of the characteristic direction $T$.\\
We explicitly note that the trajectories of the characteristic direction belong to a $(n+1)$-dimensional torus $\mathbb{T}^{n+1}$ (eventually degenerate) identified by
\begin{equation}\label{torus}
\mathbb{T}^{n+1}=\mathbb{S}^{1} \times \ldots \times \mathbb{S}^{1}=
\{z\in \Omega : |z_1|=c_1\geq 0,\ldots,|z_{n+1}|=c_{n+1}\geq 0\}
\end{equation}
and this is a particular case of the wellknown Liouville-Arnold Theorem \cite{arnold}.\\
In other words we have a symplectic toric action group on $\mathbb{C}^{n+1}$ with a fixed point at the origin.\\
Let us now consider the following explicit formula to compute the  $j$-th Levi curvature of $M$ in term of a defining
function $f$ (see \cite{ML}):
\begin{equation}\label{k1}
  L^j=-\frac{1}{\begin{pmatrix}
 n\\
 j
\end{pmatrix}}\frac{1}{|\partial f|^{j+2}}\sum_{1\leq i_1<\dots<i_{j+1}\leq n+1}
\Delta_{(i_1,\cdots,i_{j+1})}(f)
\end{equation}
for all $j=1,\dots,n$, where
\begin{equation}\label{k2}
\Delta_{(i_1,\cdots,i_{j+1})}(f)=\det\left(
\begin{array}{llll}
 0 & f_{\bar i_1} & \ldots &  f_{\overline{i}_{j+1}}\\
 f_{i_1}  & f_{i_1, \bar i_1} & \ldots & f_{i_1, \overline{i}_{j+1}} \\
  \vdots & \vdots & \ddots & \vdots \\
 f_{i_{j+1}}  & f_{i_{j+1}, \bar i_1} & \ldots & f_{i_{j+1}, \overline{i}_{j+1}}
\end{array}%
\right)
\end{equation}
If $f(z)=g(r)$ depends only on the radii $r=(r_1,\ldots,r_{n+1})$ then by a direct computation we have that $\Delta_{(i_1,\cdots,i_{j+1})}(g)$ depends only on $(r_{i_1},\cdots,r_{i_{j+1}})$. Thus all the $j$-th Levi curvatures are conserved quantities on every fixed $(n+1)$-dimensional torus $\mathbb{T}^{n+1}$: in particular they are constant along the trajectories of the characteristic direction $T$.\\
Moreover by the formula (\ref{reinhardtcharacteristiccurvature}) also the characteristic curvature $h(T,T)$ is constant on every fixed $(n+1)$-dimensional torus.
We explicitly recall that $h(T,T)$ (and all the conserved quantities as well) is constant along the trajectories of the characteristic direction $T$ but the value of the constant changes accordingly to the initial condition of the equation (\ref{eq:hamrein}).\\
Then our main result Theorem (\ref{teoreinhardt}) states that if the value of the constant $h(T,T)$ is the same on all the trajectories of the characteristic direction $T$ then $M$ is a sphere.


\addcontentsline{toc}{section}{References}



\begin{thebibliography}{100}

\bibitem{alex} {\sc A.D. Alexandrov},
\emph{A characteristic property of spheres}, Ann. Mat. Pura Appl.,
(4) 58, pag. 303-315, 1962

\bibitem{arnold} {\sc V.I. Arnold},
\emph{Mathematical Methods of Classical Mechanics}, Springer-Verlag (1989)


\bibitem{bog} {\sc A. Bogges},
\emph{CR Manifolds and the Tangential Cauchy-Riemann Complex},
Studies in Advanced Mathematics, 1991


\bibitem{hormandercomplex} {\sc L. H\"{o}rmander},   {\em An introduction to complex analysis in several variables}, North-Holland  (1973)

\bibitem{HL} {\sc J.G. Hounie, E. Lanconelli}, {\em An Alexander type Theorem for Reinhardt domains of $\mathbb{C}^2$. Recent progress on some problems in several complex variables and partial differential equations}, Contemp. Math., Amer. Math. Soc., Providence, RI,   400,129--146, 2006.

\bibitem{HL1} {\sc J.G. Hounie, E. Lanconelli}, {\em A sphere theorem for a class of Reinhardt domains with constant Levi curvature}, Forum Mathematicum  2008  20:4  , 571-586

\bibitem{jarpfl} {\sc M. Jarnicki, P. Pflug},
\emph{ First steps in several complex variables: Reinhardt domains}, European Mathematical Society, 2008


\bibitem{kling} {\sc W.Klingenberg},
\emph{Real hypersurfaces in K\"{a}hler manifolds}, Asian J. Math. 5,
no. 1, pag. 1-17, 2001


\bibitem{ML} {\sc E. Lanconelli, A. Montanari},  \emph{Pseudoconvex Fully Nonlinear Partial Differential Operators. Strong Comparison Theorems},  J. Differential Equations 202 (2004), no. 2,  306-331

\bibitem{marmon} {\sc V. Martino, A. Montanari}, \emph{Integral formulas for a class of curvature PDE's and applications to isoperimetric inequalities and to symmetry problems}, Forum Mathematicum,  vol. 22/2 (2010); p. 255 - 267

\bibitem{MaM} {\sc V. Martino, A. Montanari}, {\em On the characteristic direction of real hypersurfaces in $\mathbb{C}^{n+1}$ and a symmetry result}, Advances in Geometry, Vol. 10/3, 2010, 371 - 377

\bibitem{io} {\sc V. Martino},
\emph{La forma di Levi per ipersuperfici in
$\mathbb{C}^{N+1}$ e l'equazione di pseudocurvatura media per
grafici reali}, PhD Thesis

\bibitem{momo} {\sc R. Monti and D. Morbidelli},
\emph{Levi umbilical surfaces in complex space}, J. Reine Angew.
Math. Math. 603 (2007) 113-131

\bibitem{RE2}{\sc R. C. Reilly,} {\em Applications of the Hessian operator in a Riemann manifold} Indiana Univ. Math. J.,
{\bf 26} (1977) 459-472.

\bibitem{R}{\sc R. C. Reilly,} {\em Mean curvature, the Laplacian, and Soap Bubbles}
Amer. Math. Monthly 89 (1982), no. 3, 180--188, 197--198.

\bibitem{R2}{\sc R. C. Reilly,} {\em On the Hessian of a function and the curvatures of its graph.}
Michigan Math. J. 20 (1973), 373--383.

\end{thebibliography}
\end{document}